\documentclass[12pt]{article}

\newtheorem{theorem}{Theorem}[section]

\newtheorem{proposition}[theorem]{Proposition}

\newtheorem{remark}[theorem]{\sc Remark}

\newfam\msbfam
\font\tenmsb=msbm10  scaled \magstep1 \textfont\msbfam=\tenmsb
\font\sevenmsb=msbm7 scaled \magstep1 \scriptfont\msbfam=\sevenmsb
\font\fivemsb=msbm5  scaled \magstep1 \scriptscriptfont\msbfam=\fivemsb
\def\Bbb{\fam\msbfam \tenmsb}

\def\RR{{\Bbb R}}
\def\CC{{\Bbb C}}
\def\QQ{{\Bbb Q}}
\def\NN{{\Bbb N}}
\def\ZZ{{\Bbb Z}}
\def\II{{\Bbb I}}
\def\TT{{\Bbb T}}
\def\BB{{\Bbb B}}

\def\ss{\subseteq}
\def\ra{\rightarrow}

\def\detjac{\hbox{det}\, \hbox{Jac}_\CC \, }

 \def\HollowBox #1#2{{\dimen0=#1 \advance\dimen0 by -#2       
       \dimen1=#1 \advance\dimen1 by #2                       
        \vrule height #1 depth #2 width #2                    
        \vrule height 0pt depth #2 width #1                   
        \llap{\vrule height #1 depth -\dimen0 width \dimen1}%
       \hskip -#2                                             
       \vrule height #1 depth #2 width #2}}                   
 \def\BoxOpTwo{\mathord{\HollowBox{6pt}{.4pt}}\;}             

\def\endpf{\hfill $\BoxOpTwo$ \smallskip \\ }

\newfam\msbbfam
\font\tenmsbb=msbm10  scaled \magstep1 \textfont\msbbfam=\tenmsbb
\font\sevenmsbb=msbm7  scaled \magstep1 \scriptfont\msbbfam=\sevenmsbb
\font\fivemsbb=msbm5    scaled \magstep1 \scriptscriptfont\msbbfam=\fivemsbb

\usepackage{graphicx}

\usepackage{amsmath}

\begin{document}

\begin{center}
\Large \bf On the Bergman Projection and the Lu Qi-Keng Conjecture\footnote{{\bf Subject 
Classification Numbers:}   32A36, 32A40, 32A70, 30H20.}
\end{center}
\vspace*{.12in}

\begin{center}
Steven G. Krantz
\end{center}

\date{\today}

\begin{quote}
{\bf Abstract:}   On a reasonable class of domains in $\CC^n$, we characterize 
those holomorphic functions which continue analytically
past the boundary.  Then we give some applications of this result to holomorphic
mappings.  In addition, some new results about the Lu Qi-Keng conjecture 
are treated.
\end{quote}
\vspace*{.25in}

\markboth{STEVEN G. KRANTZ}{THE BERGMAN PROJECTION}

\section{Introduction}

Work of S. R. Bell (see, e.g., [BEL1]) has demonstrated
the importance of the Bergman kernel $K$ and Bergman projection $P$
in understanding holomorphic mappings.  In particular, Bell's
Condition $R$ has played a central role for many years.

In the present paper we give a characterization of those functions $\varphi$
on a domain $\Omega$ such that $\varphi$ continues analytically past
the boundary.  Then we give some applications of this result.

It is a pleasure to thank Harold Boas, Siqi Fu, and Emil
Straube for many helpful remarks and suggestions.

\section{Principal Results}

In what follows a {\it domain} in $\CC^n$ is a connected, open set.
Now our main result is this:

\begin{theorem}  \sl
Let $\Omega \ss \CC^n$ be a bounded, pseudoconvex domain with
real analytic boundary.  Assume that the $\overline{\partial}$-Neumann
problem on $\Omega$ is real analytic hypoelliptic.   If $\varphi$ is a holomorphic
function on $\Omega$ that continues real analytically across all boundary
points of $\Omega$, then we may find a $g \in C_c^\infty(\Omega)$ such
that $P g = \varphi$.
\end{theorem}

\begin{remark} \rm
It may be noted that, if $g \in C_c^\infty(\Omega)$, then $Pg$ automatically continuous
analytically across the boundary by the real analytic hypoellipticity of
the $\overline{\partial}$-Neumann problem and by Kohn's projection
formuls $P = I - \overline{\partial}^* N \overline{\partial}$.
In particular, since the Bergman kernel $K$ is just the Bergman
projection of the Dirac delta mass, we see that $K \, \cdot \, , \zeta)$
analytically continues across the boundary for $\zeta \in \Omega$ fixed.  
\end{remark}
\vspace*{.12in}

\noindent {\bf Proof of the Theorem:}  Let ${\cal O}(\overline{\Omega})$ denote those
functions which are holomorphic on a neighborhood of the closure
$\overline{\Omega}$ of the domain $\Omega$.  Our job is to show that
the Bergman projection $P$ maps $C_c^\infty(\Omega)$
onto ${\cal O}(\overline{\Omega})$.   This is equivalent
to showing that the adjoint mapping (which is also $P$) maps
${\cal O}^*(\overline{\Omega})$ univalently into the
dual of $C_c^\infty(\Omega)$.  The latter is of course
just the space of distributions on $\Omega$.

Now let $\lambda$ be an element of ${\cal O}^*(\overline{\Omega})$.   We need
to see that if $\lambda \ne 0$ then $P\lambda \ne 0$.  Suppose
to the contrary that
$$
\langle P \lambda , \psi \rangle = 0
$$
for every $\psi \in C_c^\infty(\Omega)$.  Then it
follows that
$$
\langle \lambda, P \psi \rangle = 0
$$
for every $\psi \in C_c^\infty(\Omega)$.  We may write this
last as
$$
\lambda \int_\Omega K(z, \zeta) \psi(\zeta) \, dV(\zeta) = 0 \, .
$$

The last displayed equation may be written as
$$
\int_\Omega \lambda_z K(z, \zeta) \psi(\zeta) \, dV(\zeta) = 0
$$
for all $\psi \in C_c^\infty(\Omega)$.   But then it
would follow that
$$
\int_\Omega \lambda_z K(z, \zeta) h(\zeta) \, dV(\zeta) = 0
$$
for every locally $L^2$ function $h$ on $\Omega$.  Hence, for each fixed $z$,
$$
\lambda_z  K(z, \, \cdot \, ) \equiv 0 \, .
$$
By earlier remarks, this is true even for $z$ in the boundary.

But this would mean that, if $b$ is any element of the Bergman space on $\Omega$,
then 
$$
\lambda b = \int_\Omega \lambda_z K(z, \zeta) b(\zeta) \, dV(\zeta) \equiv 0 \, .
$$
Hence $\lambda$ is the zero functional, which is a contradiction.  So
the adjoint of $P$ is univalent.  Hence $P$ maps
$C_c^\infty(\Omega)$ onto ${\cal O}(\overline{\Omega})$.
\endpf
\smallskip \\

\begin{remark} \rm
It would be incorrect to suppose that if $P f \in {\cal O}(\overline{\Omega})$, then
$f \in C_c^\infty(\Omega)$.  For example, if {\bf 1} denotes the function
that is identically 1 on $\Omega$ then $P {\bf 1} = {\bf 1}$.
\end{remark}

\section{An Application}

In the paper [ALE], H. Alexander proved the following striking result:

\begin{theorem} \sl
Let $\Phi$ be a proper holomorphic mapping of the unit ball $B$ in $\CC^n$, $n > 1$,
to itself.  Then in fact $\Phi$ must be a biholomorphism.
\end{theorem}

This solved a problem of longstanding, and was a conceptually important result
at the time.  It contrasts of course with the situation in $\CC^1$.  
Shortly thereafter, W. Rudin [RUD] came up with a much
more elementary proof of a more general result.   A bit later, S. Bell [BEL2]
was able to put these ideas into a more natural context and give a proof
that used key ideas from mapping theory.  He was also able to generalize the
result from the ball to a more general class of domains.

Recall now the Lu Qi-Keng conjecture (see [BOA]).  The question is whether
the Bergman kernel for a domain $\Omega \ss \CC^n$ ever vanishes on $\Omega \times \Omega$.
Thanks to work of Boas and others, the answer is known to be negative in a number
of cases.  But the answer is affirmative, for example, on a bounded, homogeneous,
complete circular domain.   A domain for which the conjecture is true
is said to have the {\it Lu Qi-Keng property}.

Here we generalize Bell's result and put his proof into a simple setting.
The main result is as follows:

\begin{theorem} \sl
Let $\Omega_1$, $\Omega_2$ be bounded, pseudoconvex domains with real
analytic boundary and each having $\overline{\partial}$-Neumann
problem that is real analytic hypoelliptic.  Also suppose that $\Omega_1$
has the Lu Qi-Keng property.  Let $\Phi: \Omega_1 \ra \Omega_2$
be a proper holomorphic mapping.  Then in fact $\Phi$ is biholomorphic.
\end{theorem}

\noindent In the proof, we shall let $P_j$ denote the Bergman projection
on $\Omega_j$.   We begin, as in the paper [BEL2], by noting three facts:
\begin{enumerate}
\item[{\bf (a)}]  If $\varphi \in C_c^\infty(\Omega_1)$ then $P_1 \varphi$ extends
to be holomorphic on a neighborhood of $\overline{\Omega_1}$.   This is immediate
from the local real analytic hypoellipticity of the $\overline{\partial}$-Neumann
operator $N$, because $P_1 = I - \overline{\partial}^* N \overline{\partial}$.
\item[{\bf (b)}]  For each monomial $z^\alpha$, there is a function $\varphi_\alpha \in C_c^\infty(\Omega_2)$
such that $P_2 \varphi_\alpha = z^\alpha$.  This is of course a direct application
of our Theorem 2.1.
\item[{\bf (c)}]  Let $u = \hbox{det}\, (\hbox{Jac}\, \Phi)$.  If $\varphi \in L^2(\Omega_2)$, then $u \cdot (\varphi \circ \Phi) \in L^2(\Omega_1)$ and
$P_1(u \cdot (\varphi \circ \Phi)) = u \cdot ((P_2 \varphi) \circ \Phi)$.   This is a standard formula of Bell, for which see
[KRA1, Ch.\ 11].
\end{enumerate}
\vspace*{.15in}

\newfam\msbfam
\font\tenmsb=msbm10  scaled \magstep1 \textfont\msbfam=\tenmsb
\font\sevenmsb=msbm7 scaled \magstep1 \scriptfont\msbfam=\sevenmsb
\font\fivemsb=msbm5  scaled \magstep1 \scriptscriptfont\msbfam=\fivemsb
\def\Bbb{\fam\msbfam \tenmsb}

\def\RR{{\Bbb R}}
\def\CC{{\Bbb C}}
\def\QQ{{\Bbb Q}}
\def\NN{{\Bbb N}}
\def\ZZ{{\Bbb Z}}
\def\II{{\Bbb I}}
\def\TT{{\Bbb T}}
\def\BB{{\Bbb B}}

\noindent {\bf Proof of the Theorem:}  Since several of the key ideas appear in [BEL2], we
merely outline the argument.

Again using Theorem 2.1 above, let $\varphi_\alpha \in C_c^\infty(\Omega_2)$ 
be such that $P_2 \varphi_\alpha = z^\alpha$.  Thus
$$
u \Phi^\alpha = u \cdot ((P_2 \varphi_\alpha) \circ \Phi ) = P_1 (u \cdot (\varphi_\alpha \circ \Phi)) \, .
$$
We note that $u \cdot (\varphi_\alpha \circ \Phi)$ is a function in $C_c^\infty(\Omega_1)$ just because
$\Phi$ is a proper mapping.  Thus Fact {\bf (a)} implies that $u \Phi^\alpha$ extends to be holomorphic in
a neighborhood of $\overline{\Omega_1}$.   Now let $z \in \partial \Omega_1$.  We have that
$u \cdot \Phi^\alpha$ belongs to the ring of germs of holomorphic functions at $z$ for all
multi-indices $\alpha$, including $\alpha = (0,0, \dots, 0)$.  Because this ring is a unique
factorization domain, we may decompose each of the functions $u \cdot \Phi^\alpha$ into a product of
powers of irreducible elements of the ring.   We take the special case $\alpha = (1, 0, 0, \dots, 0)$.  A
simple analysis of the decomposition into irreducible elements (see [BEL2]) shows that $\Phi_1$ (the first
component of $\Phi$) extends to be holomorphic in a neighborhood of $z$.  Likewise, the other components
of $\Phi$ extend to be holomorphic in a neighborhood of $z$.

Finally we must show that $\Phi$ is unbranched.  For this we use the Lu Qi-Keng hypothesis
and the standard mapping formula for the Bergman kernel.  Namely, we know
that
$$
K_1(z, \zeta) = \hbox{det} \, (\hbox{Jac}_\CC \Phi)(z) \cdot K_2(\Phi(z), \Phi(\zeta)) \cdot \hbox{det} \, (\hbox{Jac}_\CC \overline{\Phi})(\zeta) \, .
$$
Now $K_1$ does not vanish on $\Omega \times \Omega$, and a simple application of Hurwitz's theorem
allows us to conclude then that $K_1$ does not vanish on $\partial \Omega \times \Omega$ (of course $K$
is the Bergman projection of the Dirac delta mass, so it analytically continues across the boundary).
But then we can conclude that $\hbox{det} \, (\hbox{Jac}_\CC \Phi)(z)$ does not vanish.  Therefore
$\Phi$ does not branch, so it must be biholomorphic.
\endpf
\smallskip \\

\section{On the Lu Qi-Keng Conjecture}

In the paper [JPDA], D'Angelo proved the Lu Qi-Keng conjecture
for domains of the form
$$
\Omega_{1,m} = \{(z_1, z_2) \in \CC^2: |z_1|^2 + |z_2|^{2m} < 1\} \, 
$$ 
where $m$ is a positive integer.  He did so by producing an explicit
formula for the Bergman kernel.

We also note that the paper [BFS] treats domains of
the form 
$$
\{(z_1, z_2, \dots, z_n) \in \CC^n: |z_1|^{2/p_1} + |z_2|^{2/p_2} + \cdots + |z_n|^{2/p_n} < 1\}
$$
for the $p_j$ positive integers.   That paper finds domains for which the
Lu Qi-Keng conjecture fails.

It has been an open problem to decide the Lu Qi-Keng conjecture for
domains of the form
$$
\Omega_{m_1, m_2, \dots, m_n} \equiv \{(z_1, z_2, \dots, z_n) \in \CC^n: |z_1|^{2m_1} + |z_2|^{2m_2} + \cdots + |z_n|^{2m_n} < 1\} \, ,
$$
where $m_1, m_2, \dots, m_n$ are positive integers.   We do so affirmatively in the present section.

To keep notation simple, we restrict attention to dimension two.  So we concentrate on a domain
$$
\Omega_{m, n} = \{(z_1, z_2) \in \CC^2: |z_1|^{2m} + |z_2|^{2n} < 1\} \, ,
$$
for $m$, $n$ positive integers.  Seeking a contradiction, we suppose that, for $j$ a large positive integer, the domain
$$
\Omega_{jm, jn} = \{(z_1, z_2): |z_1|^{2mj} + |z_2|^{2nj} < 1\}
$$
fails the Lu Qi-Keng property.  Let $K_j$ be the Bergman kernel for
this last domain, and suppose that $K_j(z, \zeta) = 0$.  Applying a rotation
$e^{i\theta}$ in the $z_1$ variable, and using the usual transformation formula
for the Bergman kernel (see [KRA1, \S 1.4), we see that $K_j( e^{i\theta} z_1, z_2, \dots, z_n, e^{i\theta} \zeta_1, \zeta_2, \dots, \zeta_n) = 0$
for all $0 \leq \theta < 2\pi$.  Now applying the sub-mean value property for subharmonic
functions to $|K_j(e^{i\theta}z_1 , z_2, \dots, z_n, e^{i\theta}\zeta_1, \zeta_2, \dots, \zeta_n)|$, we conclude that $K_j(0, z_2, \dots, z_n, 0, \zeta_2, \dots, \zeta_n) = 0$.
We may repeat this argument in the $z_2$, $z_3$ \dots, $z_n$ variables to conclude
that $K_j(0, 0) = 0$.   

Now we notice that, as $j \ra \infty$, the domains $\Omega_{jm, jn}$ converge in the
Hausdorff metric on domains to the bidisc $D^2$.  By Ramadanov's theorem (see also [KRA2]),
the Bergman kernels on the $\Omega_{jm, jn}$ converge uniformly on compact sets to the Bergman
kernel on $D^2$.  Hence the Bergman kernel on $D^2$ has zeros.  That is a contradiction.

We conclude that, for $j$ large, the Bergman kernel for $\Omega_{jm, jn}$ has no zeros.
But now we can apply Bell's projection formula for the Bergman kernel under a proper
holomorphic covering (see, for instance, [BOA]) because $\Omega_{jm, jn}$ covers $\Omega_{j'm, j'n}$ for
$j' < j$.  And we may conclude that the Bergman kernel for $\Omega_{j'm, j'n}$ has no
zeros.  Hence the Bergman kernel for $\Omega_{m,n}$ is zero-free for any positive
integers $m$ and $n$.

It is easy to see how the proof just presented generalizes to
arbitrary $\Omega_{m_1, m_2, \dots, m_n}$ in any dimension.
		 
\section{Additional Results}

We now have the following result.

\begin{proposition} \sl
Let $\Omega_1$, $\Omega_2$ be as in Theorem 2.1.  Let $\Phi: \Omega_1 \ra \Omega_2$ be biholomorphic.
Suppose that $u$ is a function in the Bergman space
of $\Omega_2$ that analytically continues past $\partial \Omega_2$.
Then $(u \circ \Phi ) \cdot \hbox{det}\, \hbox{Jac}\, \Phi$ analytically
continues past $\partial \Omega_1$.
\end{proposition}
{\bf Proof:}  By our Theorem 2.1, there is a function $g \in C_c^\infty(\Omega_2)$
such that $u = P_2 g$.  Now we calculate:
\begin{eqnarray*}
u \circ \Phi(z) \cdot \hbox{det} \, \hbox{Jac}_\CC \, \Phi(z) & = & P_2 g \circ \Phi(z) \cdot \detjac \Phi(z) \\
                 & = & \detjac \Phi(z) \cdot \int_{\Omega_2} K_2(\Phi(z), \zeta) g(\zeta) dV(\zeta) \\
                 & = & \int_{\Omega_2} K_1(z, \Phi^{-1}(\xi)) g(\zeta) \detjac \Phi^{-1}(z)  \\
                 &   & \quad  \cdot \detjac \Phi(z) \cdot \overline{\detjac \Phi^{-1}(\xi)} \, dV(\zeta) \\
		 & = & \int_{\Omega_1} K_1(z, \xi) g(\Phi(\xi)) \detjac \Phi(\xi) \, dV(\xi) \\
		 & = & P_1 ((g \circ \Phi)(\xi) \cdot \detjac \Phi(\xi)) (z) \, .
\end{eqnarray*}
Because $\Phi$ is proper, $g \circ \Phi$ is $C_c^\infty$ hence $(g \circ \Phi) \cdot \detjac \Phi$
is $C_c^\infty$.  So we see that $u \circ \Phi(z) \cdot \detjac \Phi(z)$ is the Bergman projection on $\Omega_1$ of a 
$C_c^\infty$ function.  So it analytically continues past the boundary.
\endpf
\smallskip \\

The next result is a consequence of the proof of Theorem 2.1.

\begin{proposition} \sl
Let $\Omega_1$, $\Omega_2$ be domains as in the hypothesis of Theorem 2.1
and let $\Phi$ be a biholomorphic mapping of these domains.  Then $\Phi$
and $\Phi^{-1}$ extend analytically past the boundary of
$\Omega_1$ and $\Omega_2$ respectively.
\end{proposition}

\begin{remark} \rm
Theorems 2.1, 3.2, as well as Propositions 5.1, 5.2 apply to domains of the
form
$$
\Omega = \{ (z_1, z_2, \dots, z_n) \in \CC^n: |z_1|^{2m_1} + |z_2|^{2m_2} + \cdots + |z_n|^{2m_n} < 1\}
$$
for positive integers $m_1, m_2, \dots, m_n$.
\end{remark}

\section{Concluding Remarks}

Given any function space $X$ on a domain $\Omega$, it would be of interest
to know which functions have Bergman projection that lies in $X$.   Clearly
this set of questions is related to Bell's Condition $R$.

We hope to investigate these matters in future work.

\newpage

\noindent {\Large \sc References}
\vspace*{.2in}

\begin{enumerate}

\newfam\msbfam
\font\tenmsb=msbm10  scaled \magstep1 \textfont\msbfam=\tenmsb
\font\sevenmsb=msbm7 scaled \magstep1 \scriptfont\msbfam=\sevenmsb
\font\fivemsb=msbm5  scaled \magstep1 \scriptscriptfont\msbfam=\fivemsb
\def\Bbb{\fam\msbfam \tenmsb}

\def\RR{{\Bbb R}}
\def\CC{{\Bbb C}}
\def\QQ{{\Bbb Q}}
\def\NN{{\Bbb N}}
\def\ZZ{{\Bbb Z}}
\def\II{{\Bbb I}}
\def\TT{{\Bbb T}}
\def\BB{{\Bbb B}}

\item[{\bf [ALE]}] H. Alexander, Proper holomorphic mappings in
$\CC^n$, {\it Indiana Univ.\ Math.\ Jour.} 26(1977), 137--146.

\item[{\bf [BEL1]}] S. R. Bell, Biholomorphic mappings and the
$\overline{\partial}$-problem, {\it Ann.\ of Math.} 114(1981),
103--113.

\item[{\bf [BEL2]}] S. R. Bell, An extension of Alexander's
theorem on proper self mappings of the ball in $\CC^n$, {\it
Indiana Univ.\ Math.\ Jour.} 32(1983), 69--71.

\item[{\bf [BOA]}]  H. Boas, Lu Qi-Keng's problem, {\it J. Korean
Math.\ Soc.} 37(2000), 253--267.

\item[{\bf [BFS]}]  H. P. Boas, S. Fu, and E. J. Straube,
The Bergman kernel function:  explicit formulas and zeros,
{\it Proc.\ AMS} 127(1999), 805--811.



\item[{\bf [JPDA]}] J. P. D'Angelo, An explicit computation of
the Bergman kernel function, {\it Jour.\ Geom.\ Anal.}
4(1994), 23--34.

\item[{\bf [KRA1]}]  S. G. Krantz, {\it Function Theory of
Several Complex Variables}, 2nd. ed., American Mathematical
Society, Providence, RI, 2001.

\item[{\bf [KRA2]}] S. G. Krantz, A new proof and a
generalization of Ramadanov's theorem, {\it Complex Variables
and Elliptic Eq.} 51(2006), 1125--1128.

\item[{\bf [RUD]}] W. Rudin, Holomorphic maps that extend to be
automorphisms of a ball, {\it Proc.\ Am.\ Math.\ Soc.}
81(1981), 429--432.

   

\end{enumerate}
\vspace*{.25in}

\begin{quote}
Steven G. Krantz  \\
Department of Mathematics \\
Washington University in St. Louis \\
St.\ Louis, Missouri 63130  \\
{\tt sk@math.wustl.edu}
\end{quote}

\end{document}